\documentclass[12pt]{article}
\usepackage{amsthm,amsfonts,amssymb,amscd}

\begin{document}

\begin{center}
\large \bf Remarks on Galois rational covers
\end{center}\vspace{0.5cm}

\centerline{A.V.Pukhlikov}\vspace{0.5cm}

\parshape=1
3cm 10cm \noindent {\small \quad\quad\quad \quad\quad\quad\quad
\quad\quad\quad {\bf }\newline In this note we improve the theorem
on Galois rational covers $X\dashrightarrow V$ for primitive Fano
varieties $V$, recently proven by the author, in the two
directions: we extend to the maximum the class of Galois groups
$G$, for which the proof works, and relax the conditions that must
be satisfied by the variety $V$
--- the divisorial canonicity alone is sufficient.\vspace{0.1cm}

Bibliography: 10 items.} \vspace{1cm}

14E05, 14E07\vspace{0.1cm}

\vspace{0.3cm}

{\bf 1. Primitive Fano varieties.} All varieties considered in
this note are defined over the field of complex numbers. Recall
that a projective variety $V$ of dimension $M\geqslant 3$ is a
{\it primitive Fano variety}, if it is factorial, has at most
terminal singularities, and its anti-canonical class $(-K_V)$ is
ample and generates the Picard group, $\mathop{\rm Pic} V={\mathbb
Z} K_V$.

A primitive Fano variety $V$ is {\it divisorially canonical}, if
for every effective divisor $D\sim -nK_V$, where $n\geqslant 1$,
the pair $(V,\frac{1}{n}D)$ is canonical: for every exceptional
divisor $E$ over $V$ the inequality
$$
\mathop{\rm ord\,}\nolimits_E D\leqslant n\cdot a(E),
$$
which is opposite to the Noether-Fano inequality, holds.

The divisorial canonicity is a very strong property. For many
families of higher-dimensional Fano varieties (including, for
instance, hypersurfaces of degree $M+1$ in ${\mathbb P}^{M+1}$ for
$M\geqslant 5$) it is known that a Zariski general variety in
these families is divisorially canonical. The list of those
families is given at the end of the paper.\vspace{0.3cm}

%%%%%%%%%%%%%%%%%%%%%%%%%%%%%%%%%%%%%%%%%%%%%%%%%%%%%%%%%%%
%%%%%%%%%%%%%%%%%   section 2

{\bf 2. Galois rational covers.} Fix a divisorially canonical
primitive Fano variety $V$. A rational map
$X\stackrel{d:1}{\dashrightarrow} V$ os a finite degree
$d\geqslant 1$, where $X$ is some projective variety, is called a
{\it Galois rational cover}, is the corresponding field extension
${\mathbb C} (V)\subset {\mathbb C}(X)$ is a Galois extension. In
\cite{Pukh2020a} it was shown that if the variety $V$ (in addition
to the condition of divisorial canonicity) satisfies the following
two technical conditions:

(*1) for every anti-canonical divisor $R\in |-K_V|$, every prime
number $p\geqslant 2$ and every, possibly reducible, closed subset
$Y\subset V$ of codimension $\geqslant 2$ there is a non-singular
curve $N\subset V$, such that
$$
p\not|\, (N\cdot K_V),
$$
$N\cap Y=\emptyset$ and $N$ meets $R$ transversally at
non-singular points,

(*2) for every, possibly reducible, closed subset $Y\subset V$ of
codimension $\geqslant 2$ there is a non-singular rational curve
$N\subset V$, such that $N\cap Y=\emptyset$,

\noindent then there are no Galois rational covers
$X\stackrel{d:1}{\dashrightarrow} V$ with an abelian Galois group
$G$, $|G|=d\geqslant 2$, where $X$ is a rationally connected
variety.

If the variety $V$ is non-singular, then the condition (*2) holds
automatically, since $V$ is rationally connected, see
\cite[Chapter II]{Kol96}. For a non-singular hypersurface
$V\subset {\mathbb P}^{M+1}$ of degree $M+1$ the condition (*1) is
easy to check, see \cite[Sec. 3]{Pukh2020a}. It is not hard to
check this condition for non-singular Fano complete intersections
in the projective space, too; however, we will show below that the
conditions (*1) and (*2) are in fact unnecessary and can be
dropped.\vspace{0.3cm}

%%%%%%%%%%%%%%%%%%%%%%%%%%%%%%%%%%%%%%%%%%%%%%%%%%%%%%%%%%%
%%%%%%%%%%%%%%%%%   section 3

{\bf 3. The main result.} The aim of this note is to improve the
theorem, shown in \cite{Pukh2020a}, in the two directions:
firstly, to extend the class of Galois groups $G$, for which the
proof given in \cite{Pukh2020a} works to the maximum (answering a
question of Yu. G. Prokhorov), and, secondly, to show that the
conditions (*1) and (*2) are not needed. For a group $G$ its
commutant is denoted by the symbol $[G,G]$. If the equality
$G=[G,G]$ holds, then the group $G$ is said to be perfect. The
following claim is true.

{\bf Theorem 1.} {\it For a divisorially canonical primitive Fano
variety $V$ there are no Galois rational covers
$X\stackrel{d:1}{\dashrightarrow} V$, the Galois group $G$ of
which is not perfect ($G\neq [G,G]$), where $X$ is a rationally
connected variety.}

{\bf Proof.} We will show that the theorem stated above follows
from the proof of the main result of \cite{Pukh2020a} with minimal
additional arguments. First of all, let us consider the question,
for what class of groups the proof given in \cite{Pukh2020a}
works. Assume that $X\dashrightarrow V$ is a Galois rational cover
with some, not necessarily abelian, Galois group $G$. If $G_1\lhd
G$ is a proper non-trivial normal subgroup, then the rational map
$X\dashrightarrow V$ is the composition of the rational maps
$$
X\dashrightarrow X_1\dashrightarrow V,
$$
where ${\mathbb C} (X_1)\subset {\mathbb C}(X)$ is a Galois
extension with the Galois group $G_1$, and ${\mathbb C} (V)\subset
{\mathbb C}(X_1)$ is a Galois extension with the group $G/G_1$.
Since the image of a rationally connected variety is rationally
connected, the quotient group $G/G_1$ is abelian if and only if
$G_1\supset [G,G]$, and, finally, in every finite abelian group
there is a subgroup, the quotient group by which is a cyclic group
of a prime order, Theorem 1 is equivalent to the following claim.

{\bf Theorem 2.} {\it For a divisorially canonical primitive Fano
variety $V$ and every prime number $p\geqslant 2$ there are no
Galois rational covers $X\stackrel{p:1}{\dashrightarrow} V$, the
Galois group of which is a cyclic group of order $p$, where $X$ is
a rationally connected variety.}

Now let us show that the proof of the main result of
\cite{Pukh2020a} gives Theorem 2 without using the conditions (1*)
and (2*).\vspace{0.3cm}

%%%%%%%%%%%%%%%%%%%%%%%%%%%%%%%%%%%%%%%%%%%%%%%%%%%%%%%%%%%
%%%%%%%%%%%%%%%%%   section 4

{\bf 4. Cyclic covers of the variety $V$.} Fix a prime number
$p\geqslant 2$ and a cyclic cover $\sigma\colon X\dashrightarrow
V$ of order $p$, where $X$ is a rationally connected variety,
assuming that such covers exist. We may assume that $X$ is a
non-singular projective variety and $\sigma$ is a morphism. In
\cite[Propositions 1,2]{Pukh2020a} the following objects are
constructed:

--- a birational morphism $\varphi\colon V^+\to V$, which is a
composition of blow ups with non-singular centres, where $V^+$ is
a non-singular projective variety with the Picard group
$$
\mathop{\rm Pic} V^+={\mathbb
Z}H\oplus\mathop{\bigoplus}\limits_{i\in I} {\mathbb Z} E_i,
$$
where $H=-K_V$ is the anti-canonical class of the variety $V$, the
ample generator of the group $\mathop{\rm Pic} V$ (we omit the
pull back symbol $\varphi^*$), and $E_i$, $i\in I$, are all
$\varphi$-exceptional prime divisors on $V^+$,

--- a non-singular quasi-projective variety $U_X$, a birational
morphism $\varphi_X\colon U_X\dashrightarrow X$ and a Zariski open
subset $U\subset V^+$, such that

(i) the rational map
$$
\sigma_*=\varphi^{-1}\circ\sigma\circ\varphi_X\colon
U_X\dashrightarrow V^+
$$
extends to a morphism $\sigma_U\colon U_X\to V^+$, the image of
which is $U$,

{\rm (ii)} the inequality
$$
\mathop{\rm codim\,}((V^+\setminus U)\subset V^+)\geqslant 2
$$
holds,

{\rm (iii)} the map $\sigma_U\colon U_X\to U$ is a cyclic cover of
order $p$, branched over a non-singular hypersurface $W\subset U$.

Let $\overline{W}$ be the closure of the effective divisor $W$ in
$V^+$. Then
$$
\overline{W}\sim nH+\sum_{i\in I}\zeta_i E_i
$$
for some $n\in {\mathbb Z}_+$ and $\zeta_i\in {\mathbb Z}$. In
\cite[Sec. 3]{Pukh2020a} it was shown (and this is the key step),
that
$$
n\in \{0,1\}.
$$
It is in order to exclude these two options that the conditions
(*2) (if $n=0$) and (*1) (if $n=1$) were needed. However, we will
show that these two cases are easily excluded by means of the
explicit constructions in \cite[Sec. 5]{Pukh2020a}. This would
complete the proof of Theorem 2, which implies Theorem
1.\vspace{0.3cm}

%%%%%%%%%%%%%%%%%%%%%%%%%%%%%%%%%%%%%%%%%%%%%%%%%%%%%%%%%%%
%%%%%%%%%%%%%%%%%   section 5

{\bf 5. The explicit construction of a cyclic cover.} Since the
variety $V^+$ is obtained from $V$ by means of a sequence of blow
ups, and the codimension of the complement $V^+\setminus U$ is at
least 2, there is an open subset $U^+\subset U$, the image
$U_V=\varphi(U^+)\subset V$ of which on $V$ is an open subset, and
moreover,
$$
\mathop{\rm codim} ((V\setminus U_V)\subset V)\geqslant 2
$$
and the map $\varphi|_{U^+}\colon U^+\to U_V$ is an isomorphism.
In order to construct the subset $U^+$, one should simply remove
from $U$ all closed subsets $E_i\cap U$, $i\in I$: their image on
$V$ is of codimension $\geqslant 2$. Set
$$
U^+_X=\sigma^{-1}_U(U^+),
$$
so that $\sigma_U^+\colon U^+_X\to U^+$ (where $\sigma^+_U$ is
obviously the restriction of the morphism $\sigma_U$ onto $U^+_X$)
is the cyclic cover of order $p$, branched over a non-singular
hypersurface $W^+\subset U^+$. Identifying $U^+$ and $U_V$, we can
assume that $U^+$ is an open subset of the original variety $V$.
Obviously,
$$
U^+\subset V\setminus \mathop{\rm Sing} V.
$$
Let $\overline{W^+}\subset V$ be the closure of the effective
reduced divisor $W^+$ in $V$. We have:
$$
\overline{W^+}\sim nH.
$$
For a hypersurface in the projective space the options $n\in
\{0,1\}$ can be excluded from the purely topological grounds, but
we will give an algebro-geometric proof, using only the general
properties of the variety $V$, based on the explicit construction
of a cyclic cover given in \cite[Sec. 5]{Pukh2020a}. Let us recall
that construction. Since we are interested only what happens over
an open subset $U^+=U_V\subset V$ with a small complement in $V$,
we no longer need to consider the variety $V^+$.

Arguing as in \cite[Sec. 5]{Pukh2020a}, we construct the variety
$X_0\subset V\times {\mathbb P}^1_{(x_0:x_1)}$, given by the
equation
\begin{equation}\label{01.03.2021}
a_1x_1^p-a_0x_0^p=0,
\end{equation}
where $a_0,a_1\in H^0(V,{\cal O}_V(N))$ are sections without a
common divisor of zeros on $V$. There is a commutative diagram of
maps
$$
\begin{array}{rcccl}
  & X & \stackrel{\beta}{\dashrightarrow} & X_0 & \\
  \sigma\!\!\! & \downarrow &  & \downarrow & \!\!\!\pi \\
  & V & = & V,  &
\end{array}
$$
where the upper horizontal arrow $\beta$ is a birational map and
$\pi$ is induced by the projection of the direct product $V\times
{\mathbb P}^1$ onto the first factor. Removing from $U^+$ suitable
subsets of codimension $\geqslant 2$, we may assume that the
sections $a_0$, $a_1$ have no common zeros on $U^+$, and the
hypersurfaces
$$
\{a_0|_{U^+}=0\}\quad\mbox{and}\quad \{a_1|_{U^+}=0\}
$$
(in the set-theoretic sense) are non-singular
--- although possibly reducible. Let ${\cal T}_V$ be the set of
all prime divisors on $V$, on which one of the sections $a_0, a_1$
vanishes, so that
$$
\mathop{\bigcup}\limits_{T\in {\cal T}_V} (T\cap U^+)
$$
is a non-singular (possibly reducible) hypersurface. Set for $T\in
{\cal T}_V$
$$
\mu(T)=\mathop{\rm max} \{\mathop{\rm ord}\nolimits_T a_0,
\mathop{\rm ord}\nolimits_T a_1\}.
$$
(Precisely one of the two integers in the right hand side is
positive.) Let us show that the constructions of \cite[Sec.
5]{Pukh2020a} imply the following fact.

{\bf Proposition 1.} {\it The branch hypersurface $\overline{W^+}$
contains a divisor $T\in {\cal T}_V$ if and only if}
$p\not\mathbin{|}\, \mu(T)$.

(Over the complement
$$
U^+\setminus \mathop{\bigcup}\limits_{T\in {\cal T}_V}
(T\cap U^+)
$$
the projection $\pi$ is not ramified, and the variety $X_0$ is
non-singular; this is obvious from the equation
(\ref{01.03.2021}).)\vspace{0.3cm}

%%%%%%%%%%%%%%%%%%%%%%%%%%%%%%%%%%%%%%%%%%%%%%%%%%%%%%%%%%%
%%%%%%%%%%%%%%%%%   section 6

{\bf 6. Local modifications.} Let us prove Proposition 1,
repeating the arguments of \cite[Sec. 5]{Pukh2020a} for the open
set $U^+\subset V$ (from which we can, if necessary, remove closed
subsets of codimension $\geqslant 2$). Set ${\cal X}_1=U^+\times
{\mathbb P}^1$ and $X_1=X_0\cap \mathop{\rm pr}\nolimits_V^{-1}
(U^+)$. Let us construct a sequence of locally-trivial ${\mathbb
P}^1$-bundles over $U^+$
$$
{\cal X}_1\stackrel{\beta_1}{\leftarrow} {\cal
X}_2\stackrel{\beta_2}{\leftarrow}\cdots
\stackrel{\beta_{k-1}}{\leftarrow}{\cal X}_k,
$$
with projections $\pi_i\colon {\cal X}_i\to U^+$, in the following
way. With respect to some trivialization of the ${\mathbb
P}^1$-bundle ${\cal X}_i/U^+$ over an open set, intersecting the
divisor $T\in {\cal T}_V$, the hypersurface $X_i$
--- the strict transform of $X_1$ on ${\cal X}_i$ --- is defined
by the equation
$$
a_{i,1}x_1^p-a_{i,0}x_0^p=0,
$$
where $(x_0:x_1)$ are homogeneous coordinates on ${\mathbb P}^1$
and one of the regular functions, say $a_{i,1}$, does not vanish
on $T$. Assume that
$$
\mathop{\rm ord}\nolimits_T a_{i,0}\geqslant p.
$$
Then the birational transformation
$$
\beta_i\colon {\cal X}_{i+1}\to {\cal X}_i
$$
is the composition of the blow up of the subvariety
$$
T_i=\pi_i^{-1}(T)\cap X_i
$$
and the subsequent contraction of the strict transform of the
hypersurface $\pi_i^{-1}(T)\subset {\cal X}_i$. It is easy to
check that locally in a neighborhood of the generic point of the
divisor $T$ the hypersurface $X_{i+1}\subset {\cal X}_{i+1}$ is
defined by the equation
$$
a_{i+1,1}x_1^p-a_{i+1,0}x_0^p=0,
$$
where $a_{i+1,1}|_T\not\equiv 0$ and
$$
\mathop{\rm ord}\nolimits_T a_{i+1,0}= \mathop{\rm ord}\nolimits_T
a_{i,0} - p,
$$
see \cite[Sec. 5]{Pukh2020a}.

Now setting
$$
\mu_i(T)=\mathop{\rm max} \{\mathop{\rm ord}\nolimits_T a_{i,0},
\mathop{\rm ord}\nolimits_T a_{i,1}\}
$$
for every $i=1,\dots, k$ and $T\in {\cal T}_V$, we get that for
every $i=1,\dots, k$ there is a precisely one divisor $T(i)\in
{\cal T}_V$, such that
$$
\mu_{i+1}(T(i))=\mu_i(T(i))-p,
$$
and $\mu_{i+1}(T)=\mu_i(T)$ for all $T\neq T(i)$. For the variety
$X_k\subset {\cal X}_k$ we have
$$
\mu_k(T)\leqslant p-1
$$
for all $T\in {\cal T}_V$, and moreover, $\mu_k(T)\equiv \mu(T)
\mathop{\rm mod} p$.

Now for each $T\in {\cal T}_V$ there are three options:

(0) $\mu_k(T)=0$, and then the hypersurface $X_k$ is not ramified
over $T$ and for that reason non-singular over $T$, so that
$T\not\subset \overline{W^+}$,

(1) $\mu_k(T)=1$, and then the hypersurface $X_k$ is ramified over
$T$ and non-singular over $T$, so that $T\subset \overline{W^+}$,

(2) $\mu_k(T)\in \{2,\dots, p-1\}$, and then the variety $X_k$ has
a cuspidal singularity of the type
$$
t^p-s^{\mu_k(T)}=0
$$
along the non-singular subvariety $\pi_k^{-1}(T)\cap X_k$, in
terms of some local coordinates $t,s$ on the plane; in that case
the normalization of the variety $X_k$ or the obvious sequence of
blow ups along non-singular subvarieties, isomorphic to $T$, gives
a variety, non-singular over $T$, that covers $U^+$ cyclically,
and this cyclic cover is ramified over $T$, so that here $T\subset
\overline{W^+}$, too.

Since $\mu_k(T)\equiv \mu(T) \mathop{\rm mod} p$, the proof of
Proposition 1 is complete. Q.E.D.\vspace{0.3cm}

%%%%%%%%%%%%%%%%%%%%%%%%%%%%%%%%%%%%%%%%%%%%%%%%%%%%%%%%%%%
%%%%%%%%%%%%%%%%%   section 7

{\bf 7. Exclusion of the cases $n=0$ and $n=1$.} Let us complete
the proof of Theorem 2. Assume that $n=0$, that is to say, the
hypersurface $\overline{W^+}$ is empty. This means that
$\mu(T)\equiv 0 \mathop{\rm mod} p$ for every $T\in {\cal T}_V$.
It follows that $p\mathbin{|}\, N$ and we can ``extract the root''
from the sections $a_0$, $a_1$ (see Sec. 5): there are sections
$$
e_0, e_1\in H^0(V, {\cal O}_V \left(N/p\right)),
$$
such that $a_0=e_0^p$ and $a_1=e_1^p$. But then the equation
(\ref{01.03.2021}) takes the form
$$
e_1^px_1^p-e_0^px_0^p=\prod^p_{i=1} (e_1x_1-\zeta^i e_0x_0)=0,
$$
where $\zeta=\exp(2\pi i/p)$, that is, the variety $X_0$ is
reducible and is a union of $p$ irreducible components, covering
$V$ birationally. This is impossible. The contradiction excludes
the case $n=0$.

Assume that $n=1$. In that case there is a unique divisor $T^*\in
{\cal T}_V$, for which $p\not\mathbin{|}\, \mu(T^*)$, and
moreover, $T^*\sim -K_V$ is a ``hyperplane section'' (the ample
generator of the Picard group) of the variety $V$. Since the
sections $a_0$, $a_1$ do not vanish simultaneously on any prime
divisor, we have, say, that $a_1|_{T^*}\equiv 0$ and
$a_0|_{T^*}\not\equiv 0$, and $p\mathbin{|}\, \mu(T)$ for every
prime divisor $T\neq T^*$, $T\in {\cal T}_V$. Therefore, we get:
$\mathop{\rm ord}\nolimits_T a_1\equiv 0\mathop{\rm mod} p$ for
$T\neq T^*$ and $\mathop{\rm ord}\nolimits_{T^*} a_1\not\equiv
0\mathop{\rm mod} p$. Since $a_1$ is a section of the sheaf ${\cal
O}_V(N)$, this implies that
$$
p\not\mathbin{|}\, N.
$$
On the other hand, $\mathop{\rm ord}\nolimits_T a_0\equiv
0\mathop{\rm mod} p$ for all $T\in {\cal T}_V$. Since $a_0$ is
also a section of the sheaf ${\cal O}_V(N)$, we get that
$$
p \mathbin{|}\, N.
$$
This contradiction excludes the case $n=1$ and completes the proof
of Theorems 2 and 1.\vspace{0.3cm}

%%%%%%%%%%%%%%%%%%%%%%%%%%%%%%%%%%%%%%%%%%%%%%%%%%%%%%%%%%%
%%%%%%%%%%%%%%%%%   section 8

{\bf 8. Divisorially canonical varieties.} To conclude, we give
the list of families of Fano varieties, for a general divisor in
which divisorial canonicity is known. In \cite{Pukh05} divisorial
canonicity is shown for Zariski general smooth hypersurfaces of
degree $M+1$ in ${\mathbb P}^{M+1}$ for $M\geqslant 5$ and (smooth
Zariski general) double covers of the projective space ${\mathbb
P}^{M}$, branched over a hypersurface of degree $2M$ for
$M\geqslant 3$. In \cite{Pukh15a} this result was improved: the
divisorial canonicity was shown for Zariski general hypersurfaces
of degree $M+1$ in ${\mathbb P}^{M+1}$, with at worst quadratic
singularities of rank $\geqslant 8$, for $M\geqslant 9$, and
moreover, hypersurfaces that do not satisfy the condition of
divisorial canonicity form a subset of codimension $\geqslant
\frac12 (M-6)(M-5)-5$ in ${\mathbb P}(H^0({\mathbb P}^{M+1}, {\cal
O}_{{\mathbb P}^{M+1}}(M+1)))$. For the double covers of the space
${\mathbb P}^{M}$ in \cite{Pukh15a} a similar improvement was
shown: for $M\geqslant 10$ the double space, branched over a
Zariski general hypersurface of degree $2M$ with at worst
quadratic singularities of rank $\geqslant 4$, is divisorially
canonical, and the branch hypersurfaces, for which the
corresponding double cover is not divisorially canonical, form a
set of codimension $\geqslant \frac12 (M-4)(M-1)$ in ${\mathbb
P}(H^0({\mathbb P}^{M}, {\cal O}_{{\mathbb P}^{M}}(2M)))$.

For a Zariski general non-singular complete intersection of type
$$
d_1\cdot d_2\cdot \dots \cdot d_k
$$
in ${\mathbb P}^{M+k}$, where $2\leqslant d_1\leqslant
d_2\leqslant \dots\leqslant d_k$ and $d_1+\dots +d_k=M+k$, where
the inequality
$$
M\geqslant 2k+3
$$
holds, the divisorial canonicity was shown in \cite{Pukh18a}.
Before that paper, in \cite{Pukh06b} and \cite{EcklPukh2016} the
divisorial canonicity was shown for smaller classes of complete
intersections of index 1. In \cite{Pukh08a} the divisorial
canonicity was shown for Zariski general smooth Fano double
hypersurfaces of index 1 and dimension $\geqslant 6$.

In \cite{Pukh21a} the divisorial canonicity was established for
Fano varieties of index 1 that are $d$-sheeted covers of ${\mathbb
P}^{M}$, under the assumption that they have at worst quadratic
singularities, the rank of which is bounded from below (the bound
depends on the dimension $M$ and the degree of the cover) and
satisfy certain additional conditions of general position, and the
varieties that are not divisorially canonical form a set, the
codimension of which is bounded from below by an integer-valued
function of the parameters $d$ and $M$, which grows as $\frac12
M^2$ when $M$ grows.

Finally, for complete intersections of type $d_1\cdot d_2$ in
${\mathbb P}^{M+2}$ the divisorial canonicity was shown for the
varieties with at worst quadratic and bi-quadratic singularities,
the rank of which is bounded from below, in \cite{Pukh2021b},
under the assumption that certain additional conditions of general
position are satisfied, and for the codimension of the set of
complete intersections that do not satisfy those conditions, an
estimate, similar to the estimates above, was obtained.

\begin{flushleft}
Department of Mathematical Sciences,\\
The University of Liverpool
\end{flushleft}

\noindent{\it pukh@liverpool.ac.uk}

\end{document}